\documentstyle[12pt]{article}
\date{}
\newcommand{\Newpar}{\hspace*{8mm}}
\newcommand{\VS}{\vskip4mm}
\font\msbm = msbm10 scaled 1200
\newcommand{\baron}[1]{\mbox{\msbm #1}}
\font\msbn = msbm9 scaled 1200
\newcommand{\markiz}[1]{\mbox{\msbn #1}}

\font\eusm = eusm10 scaled 1200
\newcommand{\lord}[1]{\mbox{\eusm #1}}

\font\eufm = eufm10 scaled 1200
\newcommand{\prince}[1]{\mbox{\eufm #1}}
\font\eufn = eufm9 scaled 1200
\newcommand{\graf}[1]{\mbox{\eufn #1}}
\begin{document}

\title{On the von Neumann Inequality for Linear Matrix Functions
of Several Variables\thanks{This work was partially supported by the INTAS 
grant 93-0322-Ext.}}
\author{D.S.~Kalyuzhniy}
\maketitle

\begin{abstract}
\noindent
      {\bf Abstract.} The theorem on the existence of three commuting
contractions on a Hilbert space and of a linear homogeneous matrix
function of three independent variables for which the generalized von
Neumann inequality fails is proved.

\end{abstract}

\Newpar
      J.~von~Neumann has shown in \cite{1} that for any contractive
linear operator
$T$
on a Hilbert space (i.~e.,
$T$
such that
$\| T \| \leq 1$)
and for any polynomial
$p(z)$
of a complex variable the following inequality holds:
$$
   \| p(T) \| \leq \max\limits_{z \in \bigtriangleup} \, | p(z) |,
$$
where
$\bigtriangleup = \{ z \in {\baron C} \, : \, | z | \leq 1 \}$
denotes the closed unit disk.

\Newpar
      T.~Ando \cite{2} generalized this inequality to the case of any
two commuting contractions
$T_1^{\,}$
and
$T_2^{\,}$
on a Hilbert space and any polynomial
$p(z_1^{\,}, z_2^{\,})$
of two independent variables:
$$
     \| p(T_1^{\,}, T_2^{\,}) \|
\leq \max\limits_{z \in \bigtriangleup^2} \, | p(z) |,
$$
(here and further for any positive integer
$N$
$\bigtriangleup^N = \{ z \in {\baron C}^N \, : \,
 |z_k^{\,}| \leq 1,                            \,
 k = 1, \ldots, N \}$
denotes the closed unit polydisk).
N.~Varopoulos \cite{3}, however, has shown that already for three
commuting contractions an analogous inequality, in general, fails,
namely, he has constructed the triple of commuting contractive linear
operators
$T_1^{\,}$,
$T_2^{\,}$,
$T_3^{\,}$
on some finite--dimensional Hilbert space and the homogeneous
polynomial of second degree
$p(z_1^{\,}, z_2^{\,}, z_3^{\,})$
of three independent variables, such that
$$
   \| p(T_1^{\,}, T_2^{\,}, T_3^{\,}) \|
>  \max\limits_{z \in \bigtriangleup^3} \,
       | p(z) |.
$$
{\bf Remark 1.} The degree of a polynomial in Varopoulos's example can
not be diminished. Indeed, let
${\bf T} = \{ T_1^{\,}, \ldots, T_N^{\,} \}$
be an
$N$--tuple
of contractions on some Hilbert space
${\cal H}$
with an inner product
$\left< \cdot \, ,  \, \cdot \right>_{\cal H}^{\,}$
and let
$l(z_1^{\,}, \ldots, z_N^{\,}) = a_0^{\,} + \sum\limits_{k=1}^N
                                                a_k^{\,} z_k^{\,}$
be an arbitrary polynomial of first degree of
$N$
independent variables, i.~e., a linear scalar function on
${\baron C}^N$.
Then
\begin{eqnarray*}
   \lefteqn{
      \| l(T_1^{\,}, \ldots, T_N^{\,}) \|
=
      \bigl\|
             a_0^{\,} \, I_{\cal H}^{\,}
          +  \sum\limits^N_{k=1}
                 a_k^{\,} T_k^{\,}
      \bigr\|
 \ \, = \ \,
     \sup\limits_{ \| x \| = \| y \| = 1} \,
         \bigl|
                 \bigl<
                      \bigl(
                              a_0^{\,} \, I_{\cal H}^{\,}
                            + \sum\limits_{k=1}^N a_k^{\,} T_k^{\,}
                      \bigr) \,
                      x, y
                 \bigr>_{\cal H}
         \bigr|
          }
\qquad\qquad\hskip3mm
\\
&=&
     \sup\limits_{ \| x \| = \| y \| = 1}   \,
         \bigl|
                a_0^{\,} \,
                \left<
                      x, y
                \right>_{\cal H}^{\,}
            +   \sum\limits_{k=1}^N
                    a_k^{\,} \,
                \left<
                      T_k^{\,} \, x, y
                \right>_{\cal H}^{\,}
         \bigr|
\ \, \leq \ \,
       \max\limits_{
                      \lambda \in \bigtriangleup, \,
                      z \in \bigtriangleup^N
                   } \,
           \bigl|
                   a_0^{\,} \, \lambda
                +  \sum\limits_{k=1}^N
                       a_k^{\,} \, z_k^{\,}
           \bigr|
\\
&=&
       \max\limits_{
                      | \zeta | = 1, \,
                      z \in \bigtriangleup^N
                   } \,
           \bigl|
                    a_0^{\,} \, \zeta
                +   \sum\limits_{k=1}^N
                        a_k^{\,} \, z_k^{\,}
           \bigr|
\ \, = \ \,
       \max\limits_{z \in \bigtriangleup^N}   \,
           \bigl|
                    a_0^{\,}
                +   \sum\limits_{k=1}^N
                        a_k^{\,} \, z_k^{\,}
           \bigr|
\ \, = \ \,
       \max\limits_{z \in \bigtriangleup^N}  \,
           | l(z) |
\end{eqnarray*}
(here we used the Cauchy--Bunyakovskiy--Schwartz inequality for the
estimate of an inner product and the maximum modulus principle for
analytic functions in the disk). It should be noted that the
commutativity of an
$N$--tuple
${\bf T}$
is unessential in this calculation.
\VS
\Newpar
      However, one may consider matrix--valued polynomials (i.~e.,
polynomials with matrix coefficients) of several independent
variables and the notion of such polynomial of several commuting
contractions \cite{4}, and then, as we will show, the degree of a
polynomial for which the corresponding generalized von Neumann
inequality fails can be diminished to one. The statement of this
problem has arisen under the consideration of multiparameter passive
scattering linear systems (the definition of such systems one may
find in \cite{5}), however, it seems that it may have an independent
interest, too.

\Newpar
      Thus, let
${\bf T} = \{ T_1^{\,}, \ldots, T_N^{\,} \}$
be an
$N$--tuple
of commuting contractions on a Hilbert space
${\cal H}$,
$P(z_1^{\,}, \ldots, z_N^{\,})$
be a matrix--valued polynomial of
$N$
independent variables, i.~e.,
\begin{equation}
      P(z) = \sum\limits_{
                             t \in {\markiz Z^N_+}, \
                             |t| \leq m
                         } \,
                  A_t^{\,} \, z^t
\qquad
      (z = (z_1^{\,}, \ldots, z_N^{\,}) \in {\baron C}^N),
    \label{1}
\end{equation}
where
${\baron Z}_+^N = \{ t \in {\baron Z}^N \, : \, t_k{\,}  \geq 0,
                     \, k = 1, \ldots, N \}$;
for
$ t \in {\baron Z}^N_+$
$| t | = \sum\limits_{k=1}^N
             t_k{\,}$,
$z^t = \prod\limits_{k=1}^N
            z_k^{t_k}$,
$A_t{\,} \in M_n{\,} ({\baron C}) = {\cal L} ({\baron C}^n)$,
so that we shall identify the algebra of
$n \times n$
matrices over
${\baron C}$
with the
$C^*$--algebra
of all linear operators on the finite--dimensional Hilbert space
${\baron C}^n$.
Define an operator
\begin{equation}
      P({\bf T}) = P(T_1^{\,}, \ldots, T_N^{\,})
                 = \sum\limits_{
                                  t \in {\markiz Z}^N_+, \
                                  | t | \leq m
                               } \,
                   A_t{\,} \otimes {\bf T}^t
                 \in
                   {\cal L} ( {\baron C^n} \otimes {\cal H})
                 \cong
                   {\cal L} ({\cal H}^n)
    \label{2}
\end{equation}
(here
${\bf T}^t = \prod\limits^N_{k=1} \,
                  T_k^{t_k}$).
One may consider
$P({\bf T})$
as an element of the
$C^*$--algebra
$       M_n^{\,} ({\baron C}) \otimes$
\linebreak
$\otimes {\cal L} ({\cal H}) \cong  M_n^{\,} ({\cal L} ({\cal H}))$
of all
$n \times n$
matrices over
${\cal L} ({\cal H})$,
i.~e., as a matrix with operator entries (see \cite{4}). We shall
also use the following
\VS
\noindent
{\bf Definition 1} \cite{6}.
       An
$N$--tuple
${\bf U} = \{ U_1^{\,}, \ldots, U_N^{\,} \}$
of commuting unitary operators on a Hilbert space
${\cal K}$
is called a \/ {\it unitary dilation} \/ of an
$N$--tuple
${\bf T} = \{ T_1^{\,}, \ldots, T_N^{\,} \}$
of commuting contractive operators on a Hilbert space
${\cal H} \subset {\cal K}$
if
$$
   \forall \, t \in {\baron Z}^N_+
\qquad
   {\bf T}^t = P_{\cal H}^{\,} \, {\bf U}^t  \, \bigl| \, {\cal H},
$$
where
$P_{\cal H}^{\,}$
is an orthogonal projection onto
${\cal H}$
in
${\cal K}$.
\VS
\noindent
{\bf Remark 2.}
      It follows from the general dilation theorem  of W.~Arveson
(see \cite{7}, p.~278) that an
$N$--tuple
${\bf T} =  \{ T_1^{\,}, \ldots, T_N^{\,} \}$
of commuting contractions on a Hilbert space allows a unitary
dilation if and only if for any matrix--valued polynomial
$P(z)$
of the form (\ref{1}) the following generalized von Neumann
inequality holds:
\begin{equation}
           \| P({\bf T}) \|
    \leq
           \max\limits_{z \in \bigtriangleup^N} \,
               \| P(z) \|,
    \label{3}
\end{equation}
where
$P({\bf T})$
is defined by equality (\ref{2}). In particular, for
$N=1$
and
$N=2$
it means (see \cite{8} and \cite{2}) that inequality (\ref{3})
always holds.
\VS
\Newpar
      For
$N=3$
we shall formulate now the main statement of the present paper.
\VS
\noindent
{\bf Theorem.}
     {\it There exist a triple
${\bf T} = \{ T_1^{\,}, T_2^{\,}, T_3^{\,} \}$
of commuting contractions on some finite--dimensional Hilbert space
${\cal H}$
and a triple
${\bf A} = \{ A_1^{\,}, A_2^{\,}, A_3^{\,} \}$
of
$n \times n$
matrices over
${\baron C}$
with some integer
$n > 1$
such that a linear homogeneous matrix function
$L(z_1^{\,}, z_2^{\,}, z_3^{\,}) =$
\linebreak
$= A_1^{\,} \, z_1^{\,} + A_2^{\,} \, z_2^{\,} + A_3^{\,} \, z_3^{\,}$
satisfies the inequality
     }
\begin{equation}
      \| L ({\bf T}) \|
   >  \max\limits_{z \in \bigtriangleup^3} \,
          \| L(z) \| .
   \label{4}
\end{equation}
\hskip12mm
     For the proof we shall use the following definitions.
\VS
\noindent
{\bf Definition 2} (e.~g., see \cite{9}).
      Let
${\cal A}$
and
${\cal B}$
be
$C^*$--algebras.
A linear map
$\varphi : {\cal A} \to {\cal B}$
is called \/ {\it positive} \/ if it transforms positive elements of
${\cal A}$
to positive elements of
${\cal B}$,
i.~e.,
$a \in {\cal A}$
and
$a \geq 0$
imply
$\varphi (a) \geq 0$.
\VS
\Newpar
      Let
${\lord S}$
be a linear subspace in a
$C^*$--algebra
${\cal A}$
(possibly, coinciding with
${\cal A}$),
${\cal B}$
be a
$C^*$--algebra,
$\varphi \, : \, {\lord S} \to {\cal B}$
be a linear map. Then for every positive integer
$n$
$M_n^{\,} ({\baron C}) \otimes {\cal A} \cong M_n^{\,} ({\cal A})$
is the
$C^*$--algebra of all
$n \times n$
matrices over
${\cal A}$
and
$M_n^{\,} \otimes {\lord S}$
is a linear subspace of this
$C^*$--algebra.
If
$i d_n^{\,}$
denotes the identity map of
$M_n^{\,} ({\baron C})$
onto itself, then
$\varphi_n^{\,} = i d_n^{\,} \otimes \varphi$
is a linear map of
$M_n^{\,} ({\baron C}) \otimes {\lord S}$
into
$M_n^{\,} ({\baron C}) \otimes {\cal B}$
which transforms matrices over
${\lord S}$
to matrices over
${\cal B}$
applying a linear map
$\varphi$
element by element.
\VS
\noindent
{\bf Definition 3} \cite{9}.
      A linear map
$\varphi \, : \, {\cal A} \to {\cal B}$
is called \/ {\it completely positive } \/ if maps
$\varphi_n^{\,}$
are positive for all integer
$n \geq 1$.
\VS
\noindent
{\bf Definition 4} \cite{4}.
         A linear map
$\varphi \, : \, {\lord S} \to {\cal B}$
is called \/ {\it completely contractive} \/ if maps
$\varphi_n^{\,}$
are contractive {\rm (}i.~e.,
$\| \varphi_n^{\,} \| \leq 1${\rm)}
for all integer
$n \geq 1$.
\VS
\noindent
{\bf Proof of the Theorem.}
         It follows from Remark~2 that a desired triple of
contractions
${\bf T}$
doesn't allow a unitary dilation. Let us take such triple of
contractions from the known S.~Parrott's example \cite{10}.
Let
${\prince X}$
be an arbitrary Hilbert space of dimension more than one
(${\prince X}$
may as well be taken finite--dimensional). We shall define
${\cal H} = {\prince X} \oplus {\prince X}$
and a triple
${\bf T} = \{ T_1^{\,}, T_2^{\,}, T_3^{\,} \}$
of contractive linear operators on a Hilbert space
${\cal H}$
as block matrices:
$$
    T_1^{\,}
=
   \left(
        \begin{array}{cc}
              0           &  0        \\
              B_1^{\,}    &  0
        \end{array}
   \right),
\qquad
    T_2^{\,}
=
   \left(
        \begin{array}{cc}
              0           &  0        \\
              B_2^{\,}    &  0
        \end{array}
   \right),
\qquad
    T_3^{\,}
=
   \left(
        \begin{array}{cc}
              0           &  0        \\
              I_{\graf X}^{\,}    &  0
        \end{array}
   \right),
$$
where
$B_1^{\,}$
and
$B_2^{\,}$
are arbitrary non-commuting unitary operators on
${\prince X}$,
$I_{\graf X}^{\,}$
is the identity operator on
${\prince X}$.
Evidently,
$T_k^{\,}$
are commuting operators because
$T_k^{\,} \, T_j^{\,} = 0$
($k, j = 1, 2, 3$).
Let
${\lord S} = {\rm span} \, \{ \Lambda_1^{\,}, \Lambda_2^{\,}, {\baron I} \}$
be a three--dimensional linear subspace of the
$C^*$--algebra
$C ({\baron T}^2)$
of all continuous complex--valued functions on the two--dimensional torus
${\baron T}^2$
\linebreak
(${\baron T}^N = \{ \lambda \in {\baron C}^N \, : \,
                 | \lambda_k^{\,} | = 1, k = 1, \ldots, N \}$
denotes the
$N$--dimensional torus) generated by the functions
$\Lambda_1^{\,} (\lambda) = \lambda_1^{\,}$,
$\Lambda_2^{\,} (\lambda) = \lambda_2^{\,}$,
${\baron I} (\lambda) = 1$.
We shall define the linear map
$\varphi \, : \, {\lord S} \to {\cal L} ({\prince X})$
by the following way: for any
$\alpha$,
$\beta$,
$\gamma \in {\baron C}$
let
$$
   \varphi (\alpha \, \Lambda_1^{\,} + \beta \, \Lambda_2^{\,}
           + \gamma \, {\baron I})
=  \alpha \, B_1^{\,} + \beta \, B_2^{\,} + \gamma \, I_{\graf X}^{\,}.
$$
We shall show that the map
$\varphi$
is contractive but it is not completely contractive. According to
Remark~1,
\begin{eqnarray*}
      \|
          \varphi (
                       \alpha \, \Lambda_1^{\,}
                    +  \beta \, \Lambda_2^{\,}
                    +  \gamma \, {\baron I}
                  )
      \|
&=&
      \|
          \alpha \, B_1^{\,}
       +  \beta  \, B_2^{\,}
       +  \gamma \, I_{\graf X}^{\,}
      \|
\ \, \leq \ \,
     \max\limits_{z \in \bigtriangleup^2} \,
          |
               \alpha \, z_1^{\,}
           +   \beta  \, z_2^{\,}
           +   \gamma
          |
\\
&=&
     \max\limits_{\lambda \in {\markiz T}^2} \,
         |
              \alpha \, \lambda_1^{\,}
          +   \beta  \, \lambda_2^{\,}
          +   \gamma
         |
\ \, = \ \,
    \|
         \alpha \, \Lambda_1^{\,}
      +  \beta  \, \Lambda_2^{\,}
      +  \gamma \, {\baron I}
    \|
\end{eqnarray*}
(we used also the maximum modulus principle for analytic functions on the
bidisk). Thus,
$\varphi$
is a contractive map. Let us suppose that
$\varphi$
is a completely contractive map also. Then, by Theorem~1.2.9 \cite{4}, it can
be extended to a completely positive map
$\widetilde \varphi \, : \, C ({\baron T}^2) \to {\cal L} ({\prince X})$,
where
$\widetilde \varphi ({\baron I}) = I_{\graf X}^{\,}$.
But such map has the form \cite{9}:
$$
    \forall \,
    f \in C({\baron T}^2)
\qquad
    \widetilde \varphi (f)
=   P_{\graf X}^{\,} \, \pi (f) \, \bigl| {\prince X},
$$
where
$\pi$
is a representation of the
$C^*$--algebra
$C({\baron T}^2)$
on some Hilbert space
${\cal Y} \supset {\prince X}$,
and
$P_{\graf X}^{\,}$
is an orthogonal projection onto
${\prince X}$
in
${\cal Y}$.
For any
$x \in {\prince X}$
we have:
$$
    \| x \|
=   \| B_k^{\,} \, x \|
=   \| \varphi (\Lambda_k^{\,}) \, x \|
=   \| \widetilde \varphi (\Lambda_k^{\,}) \, x \|
=   \| P_{\graf X}^{\,} \, \pi (\Lambda_k^{\,}) \, x \|
\qquad
    (k = 1, 2).
$$
On the other hand, since
$\Lambda_k^{\,}$
are unitary elements of the
$C^*$--algebra
$C({\baron T}^2)$,
$\pi (\Lambda_k^{\,})$
are unitary operators on
${\cal Y}$.
Therefore
$\| x \| = \| \pi (\Lambda_k^{\,}) \, x \|$.
Thus
$\| P_{\graf X}^{\,} \, \pi (\Lambda_k^{\,}) \, x \|
= \| \pi (\Lambda_k^{\,}) \, x \|$,
but this is possible only if
$P_{\graf X}^{\,} \, \pi (\Lambda_k^{\,}) \, x
= \pi (\Lambda_k^{\,}) \, x$,
i.~e.,
$\pi  (\Lambda_k^{\,}) \, x \in {\prince X}$.
Since
$x$
is an arbitrary element of
${\prince X}$,
from this we obtain that
$$
    \pi (\Lambda_k^{\,}) \, \bigl| \, {\prince X}
=   P_{\graf X}^{\,} \, \pi (\Lambda_k^{\,}) \, \bigl| \, {\prince X}
=   \widetilde \varphi (\Lambda_k^{\,})
=   \varphi (\Lambda_k^{\,})
=   B_k^{\,}
\qquad
    (k = 1, 2).
$$
Hence we have
$$
     B_1^{\,} \, B_2^{\,}
=    \pi (\Lambda_1^{\,}) \, \pi (\Lambda_2^{\,}) \, \bigl| \, {\prince X}
=    \pi (\Lambda_1^{\,} \, \Lambda_2^{\,}) \, \bigl| \, {\prince X}
=    \pi (\Lambda_2^{\,} \, \Lambda_1^{\,}) \, \bigl| \, {\prince X}
=    \pi (\Lambda_2^{\,}) \, \pi (\Lambda_1^{\,}) \, \bigl| \, {\prince X}
=    B_2^{\,} \, B_1^{\,},
$$
but this contradicts to the choice of operators
$B_1^{\,}$
and
$B_2^{\,}$.
Therefore
$\varphi$
is not a completely contractive map, i.~e., there is an integer
$n > 1$
for which
$\| \varphi_n^{\,} \| > 1$.
This, in turn, means that there is a triple
${\bf A} = \{ A_1^{\,}, A_2^{\,}, A_3^{\,} \}$
of
$n \times n$
matrices over
${\baron C}$
such that
\begin{equation}
      \|
          \varphi_n^{\,}
             (
                 A_1^{\,} \otimes \Lambda_1^{\,}
               + A_2^{\,} \otimes \Lambda_2^{\,}
               + A_3^{\,} \otimes {\baron I}
             )
      \|
>
      \|
             A_1^{\,} \otimes \Lambda_1^{\,}
           + A_2^{\,} \otimes \Lambda_2^{\,}
           + A_3^{\,} \otimes {\baron I}
      \|.
    \label{5}
\end{equation}
We shall investigate separately the left--hand side and the
right--hand side of this inequality.
\begin{eqnarray*}
\lefteqn{
      \|
         \varphi_n^{\,}
             (
                 A_1^{\,} \otimes \Lambda_1^{\,}
               + A_2^{\,} \otimes \Lambda_2^{\,}
               + A_3^{\,} \otimes {\baron I}
             )
      \|
=
      \|
            A_1^{\,} \otimes B_1^{\,}
          + A_2^{\,} \otimes B_2^{\,}
          + A_3^{\,} \otimes I_{\graf X}^{\,}
      \|
        }
\hskip25mm
\\
\vphantom{\sum}
&=&
      \|
           (I_n^{\,} \otimes P_{\{ 0 \} \oplus {\graf X}})
           (
              A_1^{\,} \otimes T_1^{\,}
           +  A_2^{\,} \otimes T_2^{\,}
           +  A_3^{\,} \otimes T_3^{\,}
           )
           \, \bigl| \,
           {\baron C}^n \otimes ({\prince X} \oplus \{ 0 \})
      \|
\\
\vphantom{\sum}
&\leq&
      \|
              A_1^{\,} \otimes T_1^{\,}
           +  A_2^{\,} \otimes T_2^{\,}
           +  A_3^{\,} \otimes T_3^{\,}
      \|,
\end{eqnarray*}
where
$I_n^{\,}$
is the identity
$n \times n$
matrix,
${\prince X} \oplus \{ 0 \}$
and
$\{ 0 \} \oplus {\prince X}$
are subspaces in
${\cal H} = {\prince X} \oplus {\prince X}$
of all such vectors
$h = x_1^{\,} \oplus x_2^{\,}$
that
$x_2^{\,} = 0$
(resp.
$x_1^{\,} = 0$),
and
$P_{ \{ 0 \} \oplus {\graf X}}^{\,}$
is the orthogonal projection in
${\cal H}$
onto the second one of this subspaces. The right--hand side of
(\ref{5})
$$
   \|
        A_1^{\,} \otimes \Lambda_1^{\,}
      + A_2^{\,} \otimes \Lambda_2^{\,}
      + A_3^{\,} \otimes {\baron I}
   \|
=
   \max\limits_{\lambda \in {\baron T}^2}   \,
     \|
       A_1^{\,} \, \lambda_1^{\,}
     + A_2^{\,} \, \lambda_2^{\,}
     + A_3^{\,}
     \|,
$$
by virtue of the isomorphism between the
$C^*$--algebra
$M_n^{\,} ({\baron C}) \otimes C({\baron T}^2)$
and the
$C^*$--algebra
$C({\baron T}^2, M_n^{\,} ({\baron C}))$
of all continuous
$n \times n$
matrix functions on
${\baron T}^2$
(see, e.~g., proposition~4.7.3 \cite{11}). But
\begin{eqnarray}
      \max\limits_{\lambda \in {\baron T}^2}  \,
          \|
             A_1^{\,} \lambda_1^{\,}
           + A_2^{\,} \lambda_2^{\,}
           + A_3^{\,}
          \|
&=&
      \max\limits_{\zeta \in {\baron T}^3}    \,
          \|
             A_1^{\,} \zeta_1^{\,}
           + A_2^{\,} \zeta_2^{\,}
           + A_3^{\,} \zeta_3^{\,}
          \|
    \nonumber
\\
&=&
      \max\limits_{z \in \bigtriangleup^3}    \,
          \|
             A_1^{\,} z_1^{\,}
           + A_2^{\,} z_2^{\,}
           + A_3^{\,} z_3^{\,}
          \|,
   \nonumber
\end{eqnarray}
according to the maximum norm principle for analytic matrix
functions of several variables (see, e.~g., \cite{12}).

\Newpar
      Thus, finally we obtain the inequality
$$
    \|
        A_1^{\,} \otimes T_1^{\,}
      + A_2^{\,} \otimes T_2^{\,}
      + A_3^{\,} \otimes T_2^{\,}
    \|
>
    \max\limits_{z \in \bigtriangleup^3}   \,
        \|
           A_1^{\,} z_1^{\,}
         + A_2^{\,} z_2^{\,}
         + A_3^{\,} z_3^{\,}
        \|,
$$
which coincides with (\ref{4}) for the linear matrix--function
$L (z_1^{\,}, z_2^{\,}, z_3^{\,})
 = A_1^{\,} z_1^{\,} + A_2^{\,} z_2^{\,} + A_3^{\,} z_3^{\,}$.
The proof of the theorem is complete.

\VS
\Newpar
      Dmitriy S.~Kalyuzhniy,

\Newpar
      Department of Higher Mathematics,

\Newpar
      Odessa State Academy of Civil Engineering and Architecture,

\Newpar
      Didrihson str.~4, Odessa, UKRAINE, 270029.

\end{document}